\documentclass[]{amsproc}
\newtheorem{thm}{Theorem}[section]

\theoremstyle{definition}%
\newtheorem{defn}[thm]{Definition}

\newtheorem{example}[thm]{Example}

\newtheorem{exercise}[thm]{Exercise}

\newtheorem{observation}[thm]{Observation}

\PassOptionsToPackage{pdfauthor={Vladimir Kisil},%
    pdftitle={Tokens in Combinatorics, Analysis, and Physics},%
    pdfsubject={Lecture Notes},%
    backref=page,%
  pdfkeywords={tokens, semigroup, special functions, quantum
    propagator, wavelets}
}{hyperref}

\usepackage[breaklinks=true,%
  bookmarks=true,%
  ]{hyperref}
\usepackage{graphics}

\providecommand{\bin}[2]{{#1 \choose #2}}
\newcommand{\Space}[3][]{\ensuremath{\mathbb{#2}^{#3}_{#1}{}}}
\newcommand{\FSpace}[3][]{\ensuremath{ #2_{#3}^{#1}{}}}
\newcommand{\modulus}[2][\relax]{\left| #2 \right|\ifx#1\relax\else_{#1}\fi}
\newcommand{\scalar}[3][\relax]{\left\langle #2,#3 
        \right\rangle\ifx#1\relax\else_{#1}\fi}

\begin{document}
\title[Tokens in Combinatorics, Analysis, and 
Physics]{Tokens:\\ An Algebraic Construction Common\\ 
  in Combinatorics, Analysis, and Physics}

\author[Vladimir V. Kisil]%
{\href{http://maths.leeds.ac.uk/~kisilv/}{Vladimir
    V. Kisil}}
\thanks{On leave from the Odessa University.}
\address{%
School of Mathematics\\
University of Leeds\\
Leeds LS2\,9JT\\
UK
}


\email{\href{mailto:kisilv@maths.leeds.ac.uk}{kisilv@maths.leeds.ac.uk}}

\urladdr{\href{http://maths.leeds.ac.uk/~kisilv/}%
{http://maths.leeds.ac.uk/\~{}kisilv/}}

\begin{abstract} We give a brief account of a construction called
\emph{tokens} here, which is significant in algebra, analysis,
combinatorics, and physics. Tokens allow to express a semigroup on one
set via a semigroup convolution on another set. Therefore tokens are
similar to intertwining operators but are more flexible.
\end{abstract}
\keywords{Semigroups, hypergroups, tokens, poset, multiplicative
  functions, polynomial sequence of binomial type, integral kernel,
  wavelets, refinement equation, special functions, quantum
  propagator, path integral, quantum computing}
\subjclass{Primary 43A20; Secondary: 05A40, 81S40}
\maketitle

\tableofcontents


\par
\hfill\parbox{0.6\textwidth}{\footnotesize Pure mathematics consists of tautologies, analogous to ``men
are men,'' but usually more complicated. \par \hfil  Bertrand Russell \emph{History of
Western Philosophy}, Chap. XVI}\par

\section{Introduction}
\label{sec:introduction}

It is a fact of our specialised world that mathematics today is cut
across by several borders: between pure and applied, continuous and
discreet---just to name few most significant ones. Yet there are many
important and vivid links which spread through those ``iron
curtains''.  

\section{Newton's Binomial Formula}
\label{sec:newt-binom-form}

Let us start from the fundamental \emph{binomial formula} (``Newton's
binomial'' according to Russian terminology):
\begin{equation}\label{eq:binom}
  (x+y)^n = \sum_{k=0}^n \bin{n}{k} x^k y^{n-k}, \qquad
  \textrm{ where } \quad \bin{n}{k}=\frac{n!}{k!(n-k)!}.
\end{equation}
It has the immediate meaning: to express a power of a sum via powers
of summands. This is of a doubtful computational benefit if applied
just to real numbers.
Besides that we could find several less obvious but not least
important implications. It is worth first to restate the formula
in a more symmetric form:
\begin{equation}\label{eq:binom/n!}
  \frac{(x+y)^n}{n!} = \sum_{k=0}^n \frac{x^k}{k!}
  \frac{y^{n-k}}{(n-k)!}.
\end{equation} There are its many different interpretations, see
e.g.~\cite{Cartier00a}, but we mention just one from each combinatorics,
analysis, and algebra.

\subsection{Combinatorics}
\label{sec:combinatorics}

  The expression $x^n/n!$ \emph{counts} a number of functions from a
  set $N$ of $n$ elements to a set $X$ of $x$ elements if we do not
  distinguish functions obtained by permutations of elements of
  $N$. (This certainly true for an integer $x$ but nothing could
  prevent a mathematician from ``generalisation''). Then
  formula~\eqref{eq:binom/n!} reads as follows (see
  Figure~\ref{fig:counting}):   

  \begin{figure}[htbp]
    \begin{center}
      \includegraphics{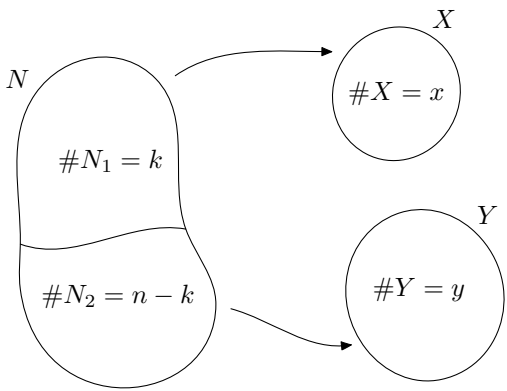}
      \caption{Counting functions $N\rightarrow X\cup Y$.}
      \label{fig:counting}
    \end{center}
  \end{figure}
  \emph{To count number of function from a set $N$ of $n$ elements to
    a union of set $X$ and $Y$ of $x$ and $y$ elements correspondingly
    split $N$ in two subsets $N'$ and $N''$ of $k$ and $n-k$ elements
    respectively, multiply number of functions from $N'$ to $X$ and
    $N''$ to $Y$ and finally sum up over all possible splittings of
    $N$.}

  The self-evidence of the above rule could serve as a proof of the
  binomial formula~\eqref{eq:binom/n!}. 

\subsection{Analysis}
\label{sec:analysis}

Let $D$ be the derivative operator or, equivalently, a linear operator
on polynomials defined by the identity:
\begin{displaymath}
  D \frac{x^n}{n!}=\frac{x^{n-1}}{(n-1)!}
  \quad \textrm{ or more generally } \quad
  D^{(k)} \frac{x^n}{n!}=\frac{x^{n-k}}{(n-k)!}.
\end{displaymath}
Consequently we could modify the formula~\eqref{eq:binom/n!} as
follows: 
\begin{eqnarray}
  \frac{(x+y)^n}{n!} 
  & = & \sum_{k=0}^n \frac{x^k}{k!} \frac{y^{n-k}}{(n-k)!} 
  \label{eq:taylorded-0}
   =  \sum_{k=0}^n \frac{x^k}{k!} 
  \left(D^{(k)}\frac{y^n}{n!}\right) \nonumber\\
  & = & \left(\sum_{k=0}^n \frac{x^k}{k!} D^{(k)}\right)
  \frac{y^n}{n!} \label{eq:taylorded-1}
   =  \left(\sum_{k=0}^\infty \frac{x^k}{k!} D^{(k)}\right)
  \frac{y^n}{n!} 
\end{eqnarray}
where firstly we extract $y^n/n!$ out of sum in \eqref{eq:taylorded-1}
because it 
is independent of $k$ and then extend summation to infinity
in~\eqref{eq:taylorded-1} because $D^{(k)} \frac{x^n}{n!}=0$ for
$k>n$.
Comparing the first~\eqref{eq:taylorded-1} and the
last~\eqref{eq:taylorded-1} lines we find that we got an expression
for the linear \emph{operator of shift} by $x$ acting on function
$y^n/n!$ in terms of \emph{linear} operator embraced
in~\eqref{eq:taylorded-1}, which is \emph{independent} of $n$. Due to
its linearity the formula is true for any linear combinations of
functions $y^n/n!$ and under suitable \emph{topological} assumptions
even for certain their limits. So we obtain the \emph{Taylor
expansion} for all functions represented by such limits:
\begin{displaymath}
  f(y+x)= \left(\sum_{k=0}^\infty \frac{x^k}{k!} D^{(k)}\right) f(y) =
  e^{xD} f(y).
\end{displaymath}

\subsection{Algebra}
\label{sec:algebra}

We use the notation $p(n,x)$ for the polynomial
$x^n/n!$. Then the formula~\eqref{eq:binom/n!} becomes:
\begin{displaymath}
  p(n,x+y)=\sum_{k=0}^n p(k,x) p(n-k,y),
\end{displaymath}
or if we define $p(n,x)\equiv 0$ for $k<0$ and any $x$ we get:
\begin{equation}\label{eq:binom-pn}
  p(n,x+y)=\sum_{k=-\infty}^\infty p(k,x) p(n-k,y).
\end{equation}
From the algebraic point of view we make the following observation. 
\begin{observation}
  \label{th:algebraic}
  The formula~\eqref{eq:binom-pn} expresses one algebraic
  operation (the \emph{sum} of two real numbers in the second argument of
  $p$'s) by another algebraic operation (a \emph{convolution} over
  integers of first arguments of $p$'s).
\end{observation}
This type of polynomials deserve a special name.
\begin{defn}\cite{BarnabBriniNicol80,Knuth99a}
  A \emph{convolution polynomials}, that is, a sequence of polynomials
  $p(n,x)$ in $x$ with $\deg p(n,x) =n$ and complex coefficients,
  satisfying the identities~\eqref{eq:binom-pn}.
\end{defn}
A cousin object is as follows:
\begin{defn}\cite{FoundationIII}
  A \emph{polynomial sequence of binomial type}, that is, a sequence
  of polynomials $p(n,x)$ in $x$ with $\deg p(n,x) =n$ and complex
  coefficients, satisfying the identities:
  \begin{equation}\label{eq:binom-type}
    p(n,x+y)=\sum_{k=-\infty}^{\infty} \bin{n}{k} p(k,x) p(n-k,y).
  \end{equation} 
\end{defn}
We prefer the explicit structure of convolution
in~\eqref{eq:binom-pn} vs.~\eqref{eq:binom-type}. The good news is
that there is an unlimited number of such polynomials.
\begin{thm}\textup{\cite{StanleyI}, \cite{Knuth99a}}
Let $f(t)$ be a formal power series:
\begin{equation}\label{eq:cumulant-gen-funct}
  f(t)=\sum_{n=1}^\infty c_n \frac{t^n}{n!},
\end{equation}
then the function $e^{xf(t)}$ is a generating function for
polynomials $p_n(x)$ of the binomial type:
\begin{equation}\label{eq:polynom-gen-func}
  \sum_{n=0}^\infty p_n(x) \frac{t^n}{n!}= e^{xf(t)}.
\end{equation} and coefficients $c_n=p_n'(0)$
in~\eqref{eq:cumulant-gen-funct} are called
\emph{cumulants}~\cite{Mattner99} for the sequence $p_n(x)$.
\end{thm}

Many examples of classic convolution polynomials obtained in this way
are listed in~\cite{RomRota78}.

\section{Semigroups and Tokens}
\label{sec:semigroups-tokens}

We would like to generalise the notion of convolution polynomials. In
order to do that we use the Observation~\ref{th:algebraic}. A proper
ground for our construction should admit an operation of convolution.
\begin{defn}
  A semigroup $C$ is called a \emph{left (right) cancellative
    semigroup} if for any $a$, $b$, $c\in C$ the identity $ca
  = cb$ ($ac=bc$) implies $a=b$. A \emph{cancellative semigroup} is
  both left and right cancellative semigroup.
\end{defn} 
We still denote the unique-if-exist solutions to equations $x \cdot
a=b$ and $a \cdot x=b$ by $[ba^{-1}]$ and $[a^{-1}b]$
correspondingly. Here the braces stress that both $[ba^{-1}]$ and
$[a^{-1}b]$ are \emph{mono}symbols and just ``$a^{-1}$'' is not
defined in general.

Of course any group is a cancellative semigroup. A simple examples
which are not groups are $\Space{N}{}$ and $\Space[+]{R}{}$ (positive
real numbers). A convolution on $\Space{N}{}$ sets was used
in~\eqref{eq:binom-pn}. There are less trivial examples.
\begin{example}\label{ex:poset}
  Let $P$ be a \emph{poset} (i.e., \emph{p}artialy \emph{o}rdered
  \emph{set}) and let $C$ denote the subset of Cartesian square
  $P\times P$, such that $(a,b)\in C$ iff $a\leq b$, $a,b\in P$. We
  can define a multiplication on $C$ by the formula:
  \begin{equation}\label{eq:interval-mult}
    (a,b)(c,d)=\left\{\begin{array}{ll}
        \textrm{undefined }, & b \not = c;\\
        (a,d), & b=c.
      \end{array} \right.
  \end{equation}
  One can see that $C$ is a c-set. If $P$ is locally
  finite, i.e., for any $a\leq b$, $a$, $b\in P$ the number of $z$ between
  $a$ and $b$ ($a\leq z\leq b$) is finite, then we can define a
  measure $d(a,b)=1$ on $C$ for any $(a,b)\in C$. With such a measure
  one defines the correct convolution on $C$:
  \begin{equation}
    h(a,b)=\int_C f(c,d)\, g([(c,d)^{-1}(a,b)])\, d(c,d)
    =\sum_{a\leq z\leq b}  f(a,z)\, g(z,b).
  \end{equation}
  The constructed algebra is the fundamental \emph{incidence algebra} in
  combinatorics~\cite{FoundationVI}.
\end{example}

Now we generalise the property of convolution
polynomials~\eqref{eq:binom-pn}.
\begin{defn}
  Let $C_1$ and $C_2$ be two c-semigroups. We will say that a function
  $t(c_1,c_2)$ on $C_1 \times C_2$ is a \emph{token}\footnote{Please
    \href{mailto:kisilv@maths.leeds.ac.uk}{tell me} if you know a
    better name for such a kind of objects.} from $C_1$ to
  $C_2$ if for any $c'_1 \in C_1$ and any $c_2, c_2' \in C_2$ we have
  \begin{equation}\label{eq:token}
    \int_{C_1} t(c_1,c_2)\, t([c_1^{-1}c_1'], c_2') \, dc_1 =
    t(c_1', c_2 c_2').
  \end{equation}
\end{defn} In fact there is another set of object, called
\emph{dissects}~\cite{Henle75}, which gives another generalisation
of~\eqref{eq:binom}, but only for a discrete sets. It is not clear at
the moment if all results about tokens~\cite{Kisil97b} could be
extended to dissects. On the other hand all examples in this paper are
described by semigroups and dissects are not necessary here.

The natural question: \emph{are tokens useful}? 
As an answer we list different examples of tokens in the next sections.

\section{Examples of Tokens}
\label{sec:examples-tokens}
Many classic object in various fields could be identified as instances
of tokens.

\subsection{Integral Kernels on Boundaries of Domains}
\label{sec:integral}
There is a clear pattern of the same structure associated with many
important integral kernels.

\begin{example}\label{ex:Poisson}
  Let $C_1$ be $\Space{R}{n}$ and $C_2$ is 
  $\Space{R}{n}\times\Space[+]{R}{}$---the ``upper half space'' in
  $\Space{R}{n+1}$. 
  For the space of harmonic function in
  $C_2=\Space{R}{n}\times\Space[+]{R}{}$  there is an integral
  representation over the boundary $C_1=\Space{R}{n}$:
\begin{displaymath}
  f(v,t)= \int_{C_1} P(u;v,t) f(u) \, du, \qquad u\in C_1, \ (v,t)\in
  C_2,\ v\in\Space{R}{n},\ t\in\Space[+]{R}{}.
\end{displaymath}
Here $P(u,v)$ is the celebrated Poisson kernel
\begin{displaymath}
  P(u;v,t)=\frac{2}{\modulus{S_n}}\frac{t}{(\modulus{u-v}^2+t^2)^{(n+1
      )/2}}
\end{displaymath}
with the property usually referred as a \emph{semigroup
  property}~\cite[Chap.~3, Prob.~1]{Akhiezer88}
\begin{displaymath}
  P(u;v+v', t+t')=\int_{C_1} P(u';v',t') P(u-u';v,t)\, du'.
\end{displaymath}
We meet the token in analysis.
\end{example}
\begin{example}\label{ex:Weierstrass}
  We preserve the meaning of $C_1$ and $C_2$ from the previous Example
  and define the Weierstrass (or Gauss-Weierstrass) kernel by the
  formula:
  \begin{displaymath}
     W(z;w,\tau)=\frac{1}{(\sqrt{2\pi \tau})^n} e^{\frac{-
        \modulus{z-w}^2}{2\tau}}=
    \frac{1}{(2\pi)^n}\int_{\Space{R}{n} }
    e^{-\frac{\tau}{2}\modulus{u}^2} e^{-(u,z-w)}\, du,
  \end{displaymath}
  where $z\in C_1$, $(w,\tau)\in C_2$. Function $W(z;w,\tau)$ is the
  fundamental solution to the heat equation~\cite[\S~2.3]{Evans98}.
  We again have~\cite[Chap.~3, Prob.~1]{Akhiezer88}
  \begin{displaymath}
    W(z;w+w',\tau +\tau')=\int_{C_1} W(z';w',\tau') W(z-z';w,\tau)\, dz.
  \end{displaymath}
  Thus we again meet a token.
\end{example} 
Two last Examples open a huge list of integral
kernels~\cite{Akhiezer88} which are tokens of
analysis.

\subsection{Multiplicative Functions on Partitions}
\label{sec:lattice-partition}

Two important examples of posets are lattices of
partitions of a set and non-crossing partitions of an ordered set.

We denote by  $\FSpace{l}{}(\Space{N}{})$ the set of sequences $f(n)$,
$n\in \Space{N}{}$, i.e. functions $\Space{N}{} \rightarrow
\Space{C}{}$. We also adopt the notion~\cite[\S~5.1]{StanleyII}
$E_f(x)$ for the \emph{exponential generating function}
\begin{displaymath}
  E_f(x)= \sum_{n=0}^\infty f(n) \frac{x^n}{n!}
\end{displaymath}
associated to $f \in \FSpace{l}{}(\Space{N}{})$.

Let us define an operation on the set of sequences
$\FSpace{l}{}(\Space{N}{})$ by means of composition corresponding
exponential generating functions:
\begin{equation}\label{eq:sequence-compos}
 h= f*g, 
 \quad \textrm{ where } \quad
 \sum_{n=0}^\infty h(n) \frac{x^n}{n!}
 =E_g(E_f(x))
 = \sum_{n=0}^\infty g(n)  \frac{(E_f(x))^n}{n!}
\end{equation}

With this operation $\FSpace{l}{}(\Space{N}{})$ is not a group: for
example $g^{-1}$ does not exist when $E_g(x)=x^2$
\cite[\S~5.4]{StanleyII}.  On the other hand
\eqref{eq:sequence-compos} makes $\FSpace{l}{}(\Space{N}{})$ a
$c$-semigroup, the uniqueness of an existing $g^{-1}f$ could be
derived similarly to proof of \cite[Prop.~5.4.1]{StanleyII}.

Let $\Phi$ be the \emph{constructor} of multiplicative
functions, i.e. for any sequence of numbers $f(n)\in
\FSpace{l}{}(\Space{N}{})$ 
and an interval $(\sigma,\pi)$ in the lattice of partition it assign a
number 
\begin{displaymath}
  \Phi(f,(\sigma,\pi))=f^{a_1}(1)f^{a_2}(2)\ldots f^{a_k}(k)\ldots, 
\end{displaymath}
where
\begin{displaymath}
  (\sigma,\pi)\simeq \Pi_1^{a_1}  \Pi_2^{a_2} \ldots \Pi_3^{a_3} \ldots.
\end{displaymath} 

We could restate the following result~\cite[Th.~5.1.11]{StanleyII}
\begin{displaymath}
  \Phi(f*g, (\sigma,\pi))=\sum_{\sigma\leq \nu\leq\pi}
  \Phi(f, (\sigma,\nu)) \Phi(g, (\nu,\pi))
\end{displaymath}
as an observation that $\Phi$ is a \emph{token} from $c$-semigroup of
poset $\Pi$ to a $c$-semigroup of $\FSpace{l}{}(\Space{N}{})$ with
the operation~\eqref{eq:sequence-compos}.

\subsection{Special Functions from Group Representations}
\label{sec:repr-theory-spec}

Let we have a representation $T$ of a group $G$ by invertible
operators in a Hilbert space $H$. 

\begin{defn}
  \label{de:matrix-elements}
  The \emph{matrix elements} $t_{jk}(g)$ of a representation $T$ of a
  group $G$ (with respect to a basis $\{{e}_j\}$ in $H$) are
  complex valued functions on $G$ defined by
  \begin{equation}\label{eq:matrix-elements}
    t_{jk}(g) = \scalar{T(g){e}_j}{{e}_k}.
  \end{equation}
\end{defn}
\begin{exercise}
  Show that \cite[\S~1.1.3]{Vilenkin68}
  \begin{enumerate}
  \item $T(g)\,{e}_k=\sum_j t_{jk}(g)\,{e}_j$.
  \item 
    \label{item:matrix-addition}
    $t_{jk}(g_1g_2)=\sum_n t_{jn}(g_1)\,t_{nk}(g_2)$.
\end{enumerate}
\end{exercise}

It is well known~\cite{Miller68,Vilenkin68} that many classic special
functions (e.g. 
\href{http://maths.leeds.ac.uk/~kisilv/courses/sp-repr.html#sec:groups-their-repr}{trigonometric
  function},
\href{http://maths.leeds.ac.uk/~kisilv/courses/sp-repr.html#ex:ker-J}{Legendre}
, Jacoby, and Hermite polynomials; Bessel, Hankel, and hypergeometric functions) appear from group representations according to the following
definition.
\begin{defn} \cite{Miller68,Vilenkin68}
  \label{de:special-function}
  A \emph{special function} associated with a
  representation $T$ of a group $G$ is a
  \hyperref[de:matrix-elements]{matrix element} $t_{ij}(g)$ of $T$.
\end{defn}
Important \emph{addition formulae}~\cite{Whittaker-Watson} for special functions are in
fact particular realisations of the simple identity:
\begin{displaymath}
  t_{jk}(g_1g_2)=\sum_n t_{jn}(g_1)\,t_{nk}(g_2).
\end{displaymath}
But this identity states that the matrix coefficients $t_{jk}(g)$
generated by a representation of a group $G$ \emph{are
tokens} from the semigroup $\Space{N}{2}$ (with a multiplication similar
to~\eqref{eq:interval-mult}) to the group $G$.

\subsection{Wavelets Refinement Equation}
\label{sec:wavelets}

An
\href{http://www.maths.leeds.ac.uk/~pmt6jrp/wavelets.formulae.ps}{orthogonal
  multiresolution of $\FSpace{L}{2}(\Space{R}{})$} (or wavelets
analysis) \cite{Daubechies92} is a chain of 
closed subspaces indexed by all integers:
\begin{displaymath}
  \cdots V_{-2} \subset V_{-1} \subset V_0 \subset V_{1} \subset V_{2}
  \cdots,
\end{displaymath}
if it has the following properties:
\begin{enumerate}
\item \emph{Completeness}. $\displaystyle\overline{\lim_{n\rightarrow
      \infty} V_n}=   \FSpace{L}{2}(\Space{R}{})$ and
  $\displaystyle\lim_{n\rightarrow -\infty} V_n =\{ 0\}$.
\item \emph{Scale Similarity}. $f_n(x) \in V_n \Leftrightarrow f_n(2x)
  \in V_{n+1}$. 
\item \emph{Translation Invariance}. $V_0$ has an orthonormal basis
  $\{\phi(x-n) \,\mid\, n\in \Space{Z}{}\}$ consisting of all integral
  translates of a single function $\phi(x)$---the \emph{mother wavelet}.
\end{enumerate}
\begin{example}
  The classic decomposition is obtained by means of the \emph{Haar
    wavelet}:
  \begin{displaymath}
    \phi(x)=\left\{ 
      \begin{array}{rl}
        1 & \textrm{ if } x \in (0,\frac{1}{2});\\
        -1 & \textrm{ if } x \in (\frac{1}{2}, 1);\\
        0 & \textrm{ otherwise}.
      \end{array}
    \right.
  \end{displaymath}
\end{example}

From the above conditions we obtain the following \emph{refinement
  equation} for the mother wavelet $\phi(x)$:
\begin{equation}\label{eq:refinement}
  \phi(x)=2\sum_{k=-\infty}^{\infty} h_k \phi(2x-k) 
  \quad \emph{ or } \quad 
  \phi\left(\frac{x}{2}\right)=\sum_{k=-\infty}^{\infty} 2h_k
  \phi(x-k).
\end{equation} 
Let us introduce a function $h(n,j)$ which is the $j$-th power
convolution over $\Space{Z}{}$ of the sequence $2h_n$ with itself. It
is easy to see~\cite[Example~3.16]{Kisil97b} that $h(n,j)$ is a token:
\begin{displaymath}
  h(n, j+j')=\sum_{k=-\infty}^{\infty} h(k,j) h(n-k,j')
\end{displaymath}
From $\Space{Z}{}$ to $\Space{N}{}$. Now
identity~\eqref{eq:refinement} just state that functions
$\psi_x(n,j)=\phi(2^jx+n)$ is a dual token (the kernel of an
associated delta family~\cite[\S~4.1]{Kisil97b}) to $h(n,j)$ for any
fixed $x\in\Space{R}{}$. Particularly we should have:
\begin{displaymath}
  \psi_x(j,n+n')=\sum_{k=-\infty}^{\infty} \psi_x(k,n) \psi_x(j-k,n'),
\end{displaymath}
or equivalent translation back to the function $\phi$:
\begin{displaymath}
  \phi(2^jx+n+n') = \sum_{k=-\infty}^{\infty}   \phi(2^kx+n)
  \phi(2^{j-k}x+n'). 
\end{displaymath}

\subsection{Quantum Mechanical Propagator}
\label{sec:quantum-mechanics}

\begin{figure}[htbp]
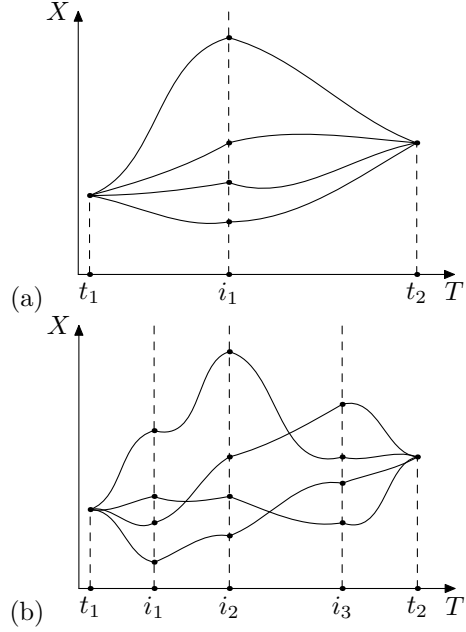

  \begin{center}
      (a)\includegraphics{tokens.2}\qquad
      (b)\includegraphics{tokens.3}
 \caption[Appearance of the path integral]
 {Repeated applications of the token property~\eqref{eq:k-proj} on
   Fig.~(a) leads to integration over paths on Fig.~(b).} 
    \label{fig:3dspectrum}
  \end{center}
\end{figure}

Let a physical system has the \emph{configuration space} $\Space{Q}{}$.
This means that we could label states of the system at any time 
$t_0 \in \Space{R}{}$ by points of $\Space{Q}{}$.  For a quantum system the
principal quantity is the \emph{propagator} \cite[\S~2.2]{FeynHibbs65},
\cite[\S~5.1]{Ryder96} $K(q_2,t_2;q_1,t_1)$---a complex valued function
defined on $\Space{Q}{} \times \Space{R}{} \times \Space{Q}{} \times
\Space{R}{}$.  It is a probability amplitude for a transition from a state
$q_1$ at time $t_1$ to $q_2$ at time $t_2$.  The \emph{probability} of this
transition is
\begin{displaymath}
  P(q_2,t_2;q_1,t_1) = \modulus{K(q_2,t_2;q_1,t_1)}^2.
\end{displaymath}
Let us fix any $t_i$, $t_1< t_i < t_2$. We assume that
\begin{equation}\label{eq:k-proj}
  K(q_2,t_2;q_1,t_1)= \int_{\Space{Q}{}} K(q_2,t_2;q_i,t_i) 
  K(q_i,t_i;q_1,t_1)\, dq_i,
\end{equation}
i.e. there exists a measure $dq$ on $\Space{Q}{}$ with the 
following property.  The system could be at time $t_i$ at any point $q_i$.
The transition amplitude $q_1 \rightarrow q_2$ is a result of all possible
transitions amplitudes $q_1 \rightarrow q_i \rightarrow q_2$ integrated
over $dq_i$.

We would like to show now that \emph{a propagator $K$ is a token in
  fact}. To do that we consider two semigroups: the semigroup
$\Space{R}{2}$ of time intervals $[t_1,t_2]$ and the semigroup
$\Space{Q}{2}$ of space intervals $[q_1,q_2]$. In the both cases the
semigroup multiplication $*$ of two intervals is given by formulas
analogous to~\eqref{eq:interval-mult}. Now we consider a propagator
$K([q_1,q_2],[t_1,t_2])=K(q_2,t_2;q_1,t_1)$ as a function on
$\Space{R}{2} \times \Space{Q}{2}$. Then the principal property of a
quantum propagator~\eqref{eq:k-proj} could be restated as
\begin{eqnarray}
  \lefteqn{K([q_1,q_2],[t_1,t']*[t',t_2])}\qquad \qquad& & \nonumber\\
  &=&\int_{\Space{Q}{2}}
  K([q_1,q'],[t_1,t']) K([q_1,q']^{-1}[q_1,q_2],[t',t_2])\,d[q_1,q'],
  \nonumber
\end{eqnarray}
which is exactly yet another realisation of the defining
property~\eqref{eq:token} of tokens.

There is an exciting but non-rigorous tool---the Feynman path
integral~\cite{FeynHibbs65}---to calculate a propagator. 
 For example it is completely strict if 
\href{http://arXiv.org/abs/math/9808040}{applied to
convolution polynomials~\cite{Kisil98d}} and expresses them via
cumulants $p_k'(0)$~\eqref{eq:cumulant-gen-funct} as
\begin{equation}\label{eq:binom-path}
  p_n(x)= n!\int\!\mathcal{D}k\mathcal{D}p\,
  \exp\!\!\int\limits_0^x
  (-ipk'+h(p))\,dt, \qquad h(p)=\sum_{k=0}^\infty p_k'(0)
  \frac{e^{ipk}}{k!}. 
\end{equation} where the first integration is taken over all possible
paths $k(t):[0,x] \rightarrow \frac{1}{x}\Space{N}{}$, such that
$k(0)=0$ and $k(x)=\frac{n}{x}$, and the path $p(t): [0,x]\rightarrow
[-{\pi},{\pi}]$ is unrestricted. It is enough to consider only paths
$k(t)$ with monotonic grow---other paths make the zero contributions.
Indeed the identity $p_l(x)\equiv 0$ for $l<0$ (made by an agreement)
implies that contribution of all paths with $k'(t)<0$ at some point
$t$ vanish. Here (and in the inner integral of~\eqref{eq:binom-path})
$k'(t)$ means the derivative of the path $k(t)$ in the distributional
sense, i.e. it is the Dirac delta function times $\frac{j}{x}$ in the
points where $k(t)$ jumps from one integer (mod $\frac{1}{x}$) value
$k(t-0)=\frac{m}{x}$ to another $k(t+0)=\frac{m+j}{x}$. Note
that~\eqref{eq:binom-path} provides an algorithm for a \emph{quantum
  computation} of a combinatorial quantity~\cite{Kisil98d}.
 
In fact the heuristic procedure defining Feynman path integral
\emph{could be applied to any token}~\cite[Rem.~4.3]{Kisil98d} with
different levels of rigour however. We hope to consider this topic
elsewhere.

\section{Conclusion}
\label{sec:some-gener-prop}

The variety of examples of tokens from different subjects given in the
previous section generates a suspicion that very little could be said
about tokens in general. Fortunately it is possible to show that many
important properties of convolution polynomials
\cite{BarnabBriniNicol80,Knuth99a} (or polynomial sequences of
binomial type \cite{FoundationIII,FoundationVI}) are
\href{http://arXiv.org/abs/funct-an/9704001/}{true for tokens in
general} \cite{Kisil97b}. Therefore all mentioned and many unmentioned
areas could benefit from the general theory of tokens on semi- and
hypergroups. We saw that tokens are
similar to intertwining operators but are more flexible. The new
frontiers for intertwining operators in functional calculus and
quantum mechanics were outlined in~\cite{Kisil02c}, it may be very
interesting to extend those ideas from intertwining operators to
tokens. 

\section*{Acknowledgments}
\label{sec:acknowledgments}

A am grateful to Franz Lehner for useful discussion of this subject
and critical remarks about this notes and my paper~\cite{Kisil97b}.

\bibliographystyle{plain}
\bibliography{abbrevmr,akisil,analyse,acombin,aphysics}

\newcommand{\noopsort}[1]{} \newcommand{\printfirst}[2]{#1}
  \newcommand{\singleletter}[1]{#1} \newcommand{\switchargs}[2]{#2#1}
  \newcommand{\irm}{\textup{I}} \newcommand{\iirm}{\textup{II}}
  \newcommand{\vrm}{\textup{V}}
  \providecommand{\cprime}{'}\providecommand{\arXiv}[1]{\eprint{http://arXiv.o%
rg/abs/#1}{arXiv:#1}}
\begin{thebibliography}{10}

\bibitem{Akhiezer88}
N.~I. Akhiezer.
\newblock {\em Lectures on Integral Transforms}.
\newblock American Mathematical Society, Providence, RI, 1988.
\newblock Translated from the Russian by H. H. McFaden.

\bibitem{BarnabBriniNicol80}
Marilena Barnabei, Andrea Brini, and Giorgio Nicoletti.
\newblock Polynomial sequences of integral type.
\newblock {\em J. Math. Anal. Appl.}, 78(2):598--617, 1980.

\bibitem{Cartier00a}
Pierre Cartier.
\newblock Mathemagics (a tribute to {L}. {E}uler and {R}. {F}eynman).
\newblock {\em S\'em. Lothar. Combin.}, 44:Art. B44d, 71 pp. (electronic),
  2000.
\newblock \MR{1814857}.

\bibitem{Daubechies92}
Ingrid Daubechies.
\newblock {\em Ten Lectures on Wavelets}, volume~61 of {\em CBMS-NSF Regional
  Conference Series in Applied Mathematics}.
\newblock Society for Industrial and Applied Mathematics (SIAM), Philadelphia,
  PA, {\noopsort{}}1992.

\bibitem{FoundationVI}
Peter Doubilet, Gian-Carlo Rota, and Richard Stanley.
\newblock On the foundations of combinatorial theory. {V}{I}. {T}he idea of
  generating function.
\newblock In Lucien~M. Le~Cam, Jerzy Neyman, and Elizabeth~L. Scott, editors,
  {\em Proceedings of the {S}ixth {B}erkeley {S}ymposium on {M}athematical
  {S}tatistics and {P}robability}, pages 267--318. University of California
  Press, Berkeley, Calif., 1972.
\newblock Reprinted in~\cite[pp.~83--134]{Rota75}
  and~\cite[pp.~148--199]{Rota95}.

\bibitem{Evans98}
Lawrence~C. Evans.
\newblock {\em Partial Differential Equations}.
\newblock American Mathematical Society, Providence, RI, 1998.

\bibitem{FeynHibbs65}
R.P. Feynman and A.R. Hibbs.
\newblock {\em Quantum Mechanics and Path Integral}.
\newblock McGraw-Hill Book Company, New York, {\noopsort{}}1965.

\bibitem{Henle75}
Michael Henle.
\newblock
  \href{http://www.ams.org/leavingmsn?url=http://links.jstor.org/sici?sici=000%
2\%2D9947\%28197502\%29202\%3C1\%3ABEOD\%3E2.0.CO\%3B2\%2DQ\%26origin=MSN}{Bin%
omial enumeration on dissects}.
\newblock {\em Trans. Amer. Math. Soc.}, 202:1--39, 1975.
\newblock \MR{50:9601}.

\bibitem{Kisil97b}
Vladimir~V. Kisil.
\newblock The umbral calculus: a model from convoloids.
\newblock {\em Z. Anal. Anwendungen}, 19(2):315--338, 2000.
\newblock \arXiv{funct-an/9704001}. \MR{2001g:05017}. \Zbl{0959.43004}.

\bibitem{Kisil02c}
Vladimir~V. Kisil.
\newblock Meeting {Descartes} and {Klein} somewhere in a noncommutative space.
\newblock In {\em Procedsings of ICMP2000}, page 25~p., 2002.
\newblock \arXiv{math-ph/0112059}.

\bibitem{Kisil98d}
Vladimir~V. Kisil.
\newblock Polynomial sequences of binomial type and path integrals.
\newblock {\em Ann. of Combinatorics}, 2002.
\newblock (Submitted)\arXiv{math/9808040}.

\bibitem{Knuth99a}
Donald~E. Knuth.
\newblock Convolution polynomials.
\newblock Preprint, Stanford Univ.,
  \href{http://www-cs-faculty.stanford.edu/~knuth/papers/cp.tex.gz}%
  {http://www-cs-faculty.stanford.edu/\~{}knuth/papers/cp.tex.gz}, February
  1996.

\bibitem{Rota95}
Joseph~P.S. Kung, editor.
\newblock {\em Gian-Carlo Rota on Combinatorics: Introductory Papers and
  Commentaries}, volume~1 of {\em Contemporary Mathematicians}.
\newblock Birkh\"auser Verlag, Boston, 1995.

\bibitem{Mattner99}
Lutz Mattner.
\newblock What are cumulants?
\newblock {\em Doc. Math.}, 4:601--622
  (\href{http://www.mathematik.uni--bielefeld.de/documenta/vol--04/18.ps.gz}{e%
lectronic}), 1999.
\newblock \MR{2001a:60017}.

\bibitem{Miller68}
Willard Miller, Jr.
\newblock {\em Lie Theory and Special Functions}.
\newblock Academic Press, New York, 1968.
\newblock Mathematics in Science and Engineering, Vol. 43.

\bibitem{FoundationIII}
Ronald Mullin and Gian-Carlo Rota.
\newblock On the foundation of combinatorial theory ({III}): Theory of binomial
  enumeration.
\newblock In B.Harris, editor, {\em Graph Theory and Its Applications}, pages
  167--213. Academic Press, Inc., New York, 1970.
\newblock Reprinted in~\cite[pp.~118--147]{Rota95}.

\bibitem{RomRota78}
S.~Roman and Gian-Carlo Rota.
\newblock The umbral calculus.
\newblock {\em Adv. in Math.}, 27:95--188, 1978.

\bibitem{Rota75}
Gian-Carlo Rota.
\newblock {\em Finite Operator Calculus}.
\newblock Academic Press, Inc., New York, 1975.

\bibitem{Ryder96}
Lewis~H. Ryder.
\newblock {\em Quantum Field Theory}.
\newblock Cambridge University Press, Cambridge, 2nd edition,
  {\noopsort{1985}}1996.

\bibitem{StanleyI}
Richard~P. Stanley.
\newblock {\em Enumerative Combinatorics. {V}ol. 1}.
\newblock Cambridge University Press, Cambridge, 1997.
\newblock With a foreword by Gian-Carlo Rota, Corrected reprint of the 1986
  original.

\bibitem{StanleyII}
Richard~P. Stanley.
\newblock {\em Enumerative Combinatorics. {V}ol. 2}.
\newblock Cambridge University Press, Cambridge, 1999.
\newblock With a foreword by Gian-Carlo Rota and Appendix 1 by Sergey Fomin.

\bibitem{Vilenkin68}
N.~Ja. Vilenkin.
\newblock {\em Special Functions and the Theory of Group Representations}.
\newblock American Mathematical Society, Providence, R. I., 1968.
\newblock Translated from the Russian by V. N. Singh. Translations of
  Mathematical Monographs, Vol. 22.

\bibitem{Whittaker-Watson}
E.~T. Whittaker and G.~N. Watson.
\newblock {\em A Course of Modern Analysis}.
\newblock Cambridge University Press, Cambridge, 1996.
\newblock An introduction to the general theory of infinite processes and of
  analytic functions; with an account of the principal transcendental
  functions, Reprint of the fourth (1927) edition.

\end{thebibliography}
\end{document}